\documentclass{amsart} 
\usepackage{amssymb} 

\numberwithin{equation}{section} 

\newtheorem{Thm}{Theorem}[section] 
\newtheorem{Prop}[Thm]{Proposition} 
\newtheorem{Lem}[Thm]{Lemma}

\theoremstyle{remark}

\theoremstyle{definition} 
\newtheorem{Def}[Thm]{Definition}

\newtheorem{Conj}[Thm]{Conjecture}

\newcommand\de{\partial}
\newcommand{\R}{\mathbb{R}}
\newcommand{\Z}{\mathbb{Z}}

\newcommand{\C}{\mathbb{C}}
  
\def\cmp#1,{{Commun.\ Math.\ Phys.\ \bf #1},}  
\def\jmp#1,{{J.\ Math.\ Phys.\ \bf #1},}  
\def\lmp#1,{{Lett.\ Math.\ Phys.\ \bf #1},}  
\def\jktr#1,{{Jour.\ of Knot Theory and its Ramification\ \bf #1},} 
\def\bams#1,{{Bull.\ Amer.\ Math.\ Soc.\ \bf #1},}
\def\agt#1,{{Algebr.\ Geom.\ Topol. \bf #1},}

\begin{document}

\title{Configuration spaces are not homotopy invariant}

\author[R. Longoni]{Riccardo~Longoni}  
\address{Dipartimento di Matematica ``G. Castelnuovo'' --- 
Universit\`a di Roma ``La  Sapienza''}   
\email{longoni@mat.uniroma1.it}  

\author[P. Salvatore]{Paolo~Salvatore}  
\address{Dipartimento di Matematica --- 
Universit\`a di Roma ``Tor Vergata''}   
\email{salvator@mat.uniroma2.it}

\begin{abstract}  
We present a counterexample to the conjecture on the 
homotopy invariance of configuration spaces. More precisely, we
consider the lens spaces $L_{7,1}$ and $L_{7,2}$, and prove that 
their configuration spaces are not homotopy equivalent by showing 
that their universal coverings have different Massey products.
\end{abstract}  

\maketitle

\section{Introduction}

The configuration space $F_n(M)$ of pairwise distinct $n$-tuples of points
in a manifold $M$ has been much studied in the literature.
Levitt reported in \cite{Levitt} as ``long-standing'' the following

\begin{Conj}
The homotopy type of $F_n(M)$, for $M$ a closed compact smooth manifold,
depends only on the homotopy type of $M$.
\end{Conj}

There was some evidence in favor: Levitt proved that the loop space 
$\Omega F_n(M)$ is a homotopy invariant of $M$. Recently Aouina and Klein 
\cite{AK} have proved that a suitable iterated suspension of $F_n(M)$ is 
a homotopy invariant.
For example the double suspension of $F_2(M)$ is a homotopy invariant.
Moreover $F_2(M)$ is a homotopy invariant when $M$ is 2-connected 
(see \cite{Levitt}). A rational homotopy theoretic version of this fact 
appears in \cite{LS}. On the other hand there is a similar situation 
suggesting that the conjecture might fail: the Euclidean configuration 
space $F_3(\R^n)$ has the homotopy type of a bundle over $S^{n-1}$ with 
fiber $S^{n-1}\vee S^{n-1}$ but it does not split as a product in general 
\cite{M2}. However the loop spaces of $F_3(\R^n)$ and of the product 
$S^{n-1} \times (S^{n-1}\vee S^{n-1})$ are homotopy equivalent and also 
the suspensions of the two spaces are homotopic.

Lens spaces provide handy examples of manifolds which are homotopy equivalent 
but not homeomorphic. The first of these examples are $L_{7,1}$ and $L_{7,2}$. 
We show that their configuration spaces $F_2(L_{7,1})$ and $F_2(L_{7,2})$ are 
not homotopy equivalent. After recalling some definition, we will describe the 
universal coverings of these configuration spaces.
Such coverings can be written as bundles with same base and fiber, but
the first splits and the second does not.
We will show that Massey products are all zero in the first case,
while there exists a nontrivial Massey product in the second case.

This means that $F_2(L_{7,1})$ is not homotopy equivalent to $F_2(L_{7,2})$. 
Finally we will extend this result by showing that  $F_n(L_{7,1})$ is not 
homotopy equivalent to $F_n(L_{7,2})$ for any $n \geq 2$.
The same result holds for unordered configuration spaces.

\section{Configuration spaces of lens spaces}

The lens spaces are 3-dimensional oriented manifolds defined as
\[
L_{m,n}:=S^3/\Z_m=\left\{(x_1,x_2)\in\C\times \C \left|\, |x_1|^2 + |x_2|^2 =1 
\right.\right\}/\Z_m
\]
where the group action is defined by 
$\zeta \left(\left(x_1,x_2\right)\right) =
(e^{2\pi i/m} x_1,e^{2\pi i n/m} x_2)$,
and $\zeta$ is the generator of $\Z_m$.
It is known (see e.g. \cite{R}) 
that $L_{7,1}$ and $L_{7,2}$ are homotopy equivalent, though not 
homeomorphic. 

For any topological space $M$, let $F_n(M)$ be the configuration space 
of $n$ pairwise distinct points in $M$, namely $F_n(M):=M^{n}\setminus
(\bigcup\Delta)$ where $\bigcup\Delta$ is the union of all diagonals. 
We first want to compute the fundamental group of $F_2(L_{7,1})$ and 
$F_2(L_{7,2})$. Observe that $S^3$ is the universal covering of $L_{7,j}$, 
for $j=1,2$, and therefore the fundamental group of $L_{7,j}$ is $\Z_7$.
Then $\pi_1(F_2(L_{7,j})) = \Z_7\times \Z_7$ because $\pi_1(L_{7,j}\times 
L_{7,j}) = \Z_7\times \Z_7$ and removing the diagonal, which is a 
codimension 3 manifold, does not change the fundamental group.

The universal coverings $\widetilde F_2(L_{7,1})$ and $\widetilde 
F_2(L_{7,2})$ are the so-called ``orbit configuration spaces'' and are 
given by pairs of points $(x,y)$ of $S^3$ which don't lie on the same 
orbit, i.e., $x\neq g(y)$ for any $g \in \Z_7$.

In the rest of the paper we identify $\Z_7$ to the group of 7th complex 
roots of unity, and we use the symbol $\zeta^t, \, t \in \R$, to denote the
complex number $e^{2\pi it/7}$. 

The first universal covering has a simple structure, namely we have the 
following
\begin{Prop}
$\widetilde F_2(L_{7,1})$ is homotopy equivalent
 to $\vee_6 S^2 \times S^3$.
\end{Prop}
\begin{proof}
It is convenient to interpret $S^3$ as the space of quaternions of unitary 
norm. Then the action of $\Z_7$ on $S^3=\widetilde{L_{7,1}}$ is the left 
translation by the subgroup $\Z_7 \subset \C\subset\mathbb H$. We define a 
map $\widetilde F_2(L_{7,1}) \to (S^3\setminus \Z_7)\times S^3$ by sending 
$(x,y)$ to $(xy^{-1},y)$. This is a homeomorphism since $x\neq \zeta^k (y) 
=\zeta^k y$ is equivalent to $xy^{-1}\neq \zeta^k$ for any 7th root of 
unity $\zeta^k, k \in \{0,\dots,6\}$.
Finally we observe that $S^3$ minus a point is $\R^3$ and hence 
$S^3\setminus\Z_7$ is homotopic to the wedge of six 2-dimensional 
spheres. 
\end{proof}

\section{Massey products}

We briefly recall the definition of Massey products for a topological 
space $X$ (see \cite{M}). Let $x,y,z\in H^*(X)$ such that $x\cup y=y
\cup z=0$. If we choose singular cochain representatives $\bar x,\bar y,\bar z
\in C^*(X)$ then we have that $ \bar x\cup \bar y = d Z$ and 
$\bar y\cup \bar z = d X$ for some cochains $Z$ and $Y$. 
Notice that $$d(Z\cup\bar z - (-1)^{\deg(x)}\bar x
\cup X) = (\bar x\cup\bar y\cup\bar z - \bar x\cup\bar y
\cup\bar z) =0,$$ and hence we can define $\langle x,y,z \rangle$ to be 
the cohomology class of $Z\cup\bar z -
(-1)^{\deg(x)}\bar x\cup X$. Since the choice of $Z$ and $X$ is not
unique, the Massey product $\langle x,y,z\rangle$ is well defined only 
in $H^*(X)/\langle y,z\rangle$ where $\langle y,z\rangle$ is the ideal 
generated by $y$ and $z$. Clearly Massey products are homotopy invariants.
A rational homotopy theoretic version of the following definition is in 
\cite{GM}.

\begin{Def}
A space $X$ is formal if the singular cochain complex
$C^*(X)$ is quasi-isomorphic to $H^*(X)$ as
augmented differential graded ring.
\end{Def}

This means there is a zig-zag of homomorphisms inducing isomorphism 
in cohomology and connecting $H^*(X)$ and $C^*(X)$. It is easy to see 
that spheres are formal. Moreover wedges and products of formal spaces 
are formal. By construction all Massey products on the cohomology of a 
formal space vanish. This in turn implies the following result

\begin{Prop}
All Massey products in the cohomology of $\widetilde F_2(L_{7,1})$ are trivial.
\end{Prop}

We deduce that in order to prove that $\widetilde F_2(L_{7,1})$ and 
$\widetilde F_2(L_{7,2})$ are not homotopy equivalent, we only need to
construct a nontrivial Massey product in the cohomology 
$\widetilde F_2(L_{7,2})$.

\section{Nontrivial Massey product for $\widetilde F_2(L_{7,2})$}

The projection onto the first coordinate gives $\widetilde F_2(L_{7,2})$
the structure of a bundle over $S^3$ with fiber 
$S^3\setminus\Z_7 \simeq \vee_6 S^2$ that admits a section. 
It follows that the Serre spectral sequence collapses and the 
cohomology ring splits as a tensor product,
so that it does not detect the nontriviality of the bundle.
In particular we have that $ H^2(\widetilde F_2(L_{7,2})) \cong \Z^6$ and
$H^4(\widetilde F_2(L_{7,2}))=0$. This in turn implies that the Massey
product of any triple in $H^2$ is well defined.

We want to compute Massey products ``geometrically'', namely using 
intersection theory on the Poincar\'e dual cycles as in \cite{M}.

Let us define the embedded ``diagonal'' 3-spheres
$\Delta_k\subset S^3\times S^3$, for $k=0,\ldots 6$, by
$\Delta_k :=\{ (x,\zeta^k(x))\,|\, x \in S^3 \}$.
Clearly $\Delta_0$ is the standard diagonal. The space $\widetilde 
F_2(L_{7,2})$ is the complement of the union of the diagonals
\[
\widetilde F_2(L_{7,2})=(S^3\times S^3)\setminus\left(\coprod_{k=0}^6\Delta_k
\right).
\]
By Poincar\'e duality we have the isomorphism
\[
H^p\left((S^3\times S^3)\setminus\left(\coprod_{k=0}^6\Delta_k
\right)\right)\cong H_{6-p}\left(S^3\times S^3, \left(\coprod_{k=0}^6\Delta_k
\right)\right).
\]
Under this identification the cup product in cohomology corresponds
to the intersection product in homology. 

We observe that there exists an isotopy $\mathcal H_k:S^3 \times [0,1] 
\to S^3 \times S^3$ (where $k$ is considered mod 7) defined by 
$\mathcal H_k((x_1,x_2),t) = ((x_1,x_2),(\zeta^{k-1+t}x_1,
\zeta^{2(k-1+t)}x_2))$. 
The images of $\mathcal H_k$ at times 0 and 1 are respectively $\Delta_{k-1}$
and $\Delta_k$, and the full image of $\mathcal H_k$ is a submanifold 
$A_k \subset S^3 \times S^3$
which represents an element in $H_4\left(S^3\times S^3, 
\left(\coprod_{k=0}^6\Delta_k \right)\right)$ Poincar\'e dual to
a class $a_k \in H^2(\widetilde F_2(L_{7,2})).$
 By using the Mayer-Vietoris
sequence one can easily see that the classes $a_k$ span 
$H^2(\widetilde F_2(L_{7,2}))$
under the relation $\sum_{k=0}^6 a_k =0$.
The main result of the paper is the following

\begin{Thm}
The Massey product $\left\langle a_4,a_1,a_2+a_6\right\rangle$
contains the class $a_2 \cup \iota$ and hence is nontrivial.
\end{Thm}

\begin{proof}
It is easy to check that $A_k$ intersects only $A_{k+3}$ and $A_{k+4}$ where
again $k$ is considered mod 7. Hence in the computation of $\left\langle 
a_4,a_1,a_2+a_6\right\rangle$ we must check the following

\begin{Lem}
\label{lemma:trasv}
The submanifolds $A_1$ and $A_4$ intersect transversally and
$$S^1 \times [0,1] \cong A_1 \cap A_4  = \left\{ \left. \left( (0,x_2),
(0,\zeta^{\lambda}x_2) \right)
\right| \, |x_2|=1,\, \lambda\in [0,1] \right\}.$$
\end{Lem}

\begin{proof}
We only need to verify that the tangent spaces to $A_1$ and $A_4$
at the point $\left((0,x_2),(0,\zeta^{\lambda}x_2)\right)$ span a six 
dimensional vector space. Recall that we are representing points in 
$S^3$ as elements $(x_1,x_2)$ in $\C\times \C$ such that 
$|x_1|^2+|x_2|^2=1$, and hence tangent vectors at $(0,x_2)$ are real linear 
combinations of the vectors $(1,0)$, $(i,0)$  
and $(0,ix_2)$. These immediately give rise to the following tangent vectors
to $A_1$ at $\left((0,x_2),(0,\zeta^{\lambda}x_2)\right)$:
$$\left((1,0),(\zeta^{\lambda/2},0)\right),\quad \left((i,0),
(i\zeta^{\lambda/2},0)\right), \quad 
\left((0,ix_2),(0,i\zeta^{\lambda}x_2)\right)$$ 
and to the following tangent vectors to $A_4$ at the same point:
$$\left((1,0),(-\zeta^{\lambda/2},0)\right),\quad \left((i,0),
(-i\zeta^{\lambda/2},0)\right),\quad 
\left((0,ix_2),(0,i\zeta^{\lambda}x_2)\right).$$
Finally consider the path in $A_1 \cap A_4$ given by 
$$s\mapsto \left((0,x_2),(0, \zeta^{\lambda + s}x_2\right).$$ Its 
derivative for $s=0$ gives, up to a scalar factor, the vector 
$\left((0,0),(0,i \zeta^{\lambda}x_2)\right)$.
By a simple inspection one sees that the linear space spanned by these 
vectors is six dimensional.
\end{proof}

Let us consider the closed 2-disc
$$D_2 = \{(r,x)\, |\, 0\le r\le 1, r^2+|x|^2=1, \, x \in \C \} \subset S^3.$$

\begin{Lem}
The intersection $A_1 \cap A_4$ is the relative boundary of the 3-manifold
\[
D_2 \times [0,1] \cong X_{14}:=
 \left\{\left.\left((r,x),(\zeta^{4t}r,\zeta^{t}x)\right)
\right|\, (r,x) \in D_2,\, 0\le t \le 1 \right\}.
\]
\end{Lem}

\begin{proof}
The pieces of the boundary of $X_{14}$ correspond to $r=0$, $t=0$ and $t=1$.
 Clearly
$\de_{r=0}X_{14} = A_1\cap A_4$. If we now show that the other pieces 
belong to one of the diagonals $\Delta_k$, the Lemma is proved. Since
$\zeta^k=\zeta^{k+7}$ we have
\begin{align*}
\de_{t=0}X_{14} & = \left\{\left((r,x),(r,x)\right)\right\}
\subset\Delta_0\\
\de_{t=1}X_{14} & = \left\{\left((r,x),(\zeta^4 r,\zeta x)\right)\right\}
\subset\Delta_4.
\end{align*}
\end{proof}

The next step is to find the intersection of $X_{14}$ with $A_2$ and $A_6$.
We observe that the inclusion $S^3 \to S^3 \times S^3$ sending $x$ to
$(1,x)$ represents the generator of $H_3\left(S^3\times S^3, 
\coprod_{k=0}^6\Delta_k \right) \cong \Z$ Poincar\'e dual to a class
$\iota \in H^3(
\widetilde F_2(L_{7,2})
 ) \cong \Z.$

\begin{Lem}
The manifolds $X_{14}$ and $A_6$ do not intersect.
Moreover $X_{14}$ and $A_2$ intersect transversally and
$X_{14}\cap A_2= A_2 \cap S^3$
is Poincar\'e dual to the class $a_2 \cup \iota$.
\end{Lem}

\begin{proof}
The intersection of $X_{14}$ with $A_6$ is given by the solution to the 
system of equations
\[
\begin{cases}
\zeta^{4t} r = \zeta^{5+s} r\\
\zeta^{t} x_2 = \zeta^{10+2s} x_2
\end{cases}
\]
for $0\le r\le 1$, $r^2+|x|^2=1$, $0\le t \le 1$ and $0\le s \le 1$. 
If we equate the exponents of the $\zeta$'s in the first and in the 
second equation we immediately see that there are no solutions for 
$0\le t \le 1$.

The intersection of $X_{14}$ with $A_2$ is given by the solution to the 
system of equations
\[
\begin{cases}
\zeta^{4t} r = \zeta^{1+s} r\\
\zeta^{t} x = \zeta^{2+2s} x
\end{cases}
\]
which has solutions $\left((1,0),(\zeta^{1+s},0)\right)$ where $0\le s
\le 1$. In fact, from the second equation we get the equation $t=2+2s$ 
(mod 7), which has no solution for $0\le t\le 1$. Therefore we must have 
$x=0$ and $r=1$. From the first equation we have that $\zeta^{4t} = 
\zeta^{1+s}$ which implies $t=(1+s)/4$. Therefore $X_{14}\cap A_2$ is a 
path connecting $\Delta_1$ with $\Delta_{2}$ which equals
$A_2\cap S^3$. 

Finally we have to check transversality for $X_{14}$ and $A_2$. By 
repeating the arguments of Lemma~\ref{lemma:trasv}, we deduce that the 
tangent space to $A_2$ at $((1,0),(\zeta^{1+s},0)) = ((1,0),(\zeta^{4t},0))$ 
is spanned by $((i,0)(i\zeta^{1+s},0))$, $((0,1),(0,\zeta^{2+2s}))$, 
$((0,i),(0,i\zeta^{2+2s}))$ and $((0,0),(i\zeta^{1+s},0))$ while the 
tangent space to $X_{14}$ at the same point
is spanned by $((0,1),(0,\zeta^{t}))$, $((0,i),(0,i\zeta^{t}))$ and 
$((0,0),(i\zeta^{4t},0))$. These vectors clearly span a six dimensional 
space.
\end{proof}

This concludes the proof since $a_2 \cup \iota$ does not belong 
to the subspace generated by $a_4 \cup \iota$ and $(a_2+a_6)\cup \iota$ 
in $$
 H^5(\widetilde F_2(L_{7,2})) =
\left\langle
a_k \cup \iota |\, k=0,\dots,6\right\rangle \left/
\sum_{k=0}^6 a_k \cup \iota \right. .
$$
\end{proof}

\section{Generalizations}

We extend our result to the $n$ points configuration space, namely we 
have that $F_n(L_{7,1})$ is not homotopic to $F_n(L_{7,2})$.
The universal covering $\widetilde{F}_n(L_{7,j})$ is the orbit configuration
space of $n$-tuples of points in $S^3$ lying in pairwise distinct 
$\Z_7$-orbits.
The forgetful map $(x_1,\dots,x_n) \mapsto (x_1,x_2)$ defines a bundle
$\widetilde{F}_n(L_{7,j}) \to \widetilde{F}_2(L_{7,j})$
 which admits a section. For example
the values $x_3,\dots,x_n$ of the section are pairwise distinct points
very close to 1 multiplied by $x_1$.  

By naturality we deduce that $\widetilde{F}_n(L_{7,2})$ has a nontrivial
Massey product on $H^2$. On the other hand right multiplication by $x_1^{-1}$
induces a product decomposition $\widetilde{F}_n(L_{7,1}) = 
S^3 \times Y_{n-1}$,
where $Y_{n-1}$ is the $n-1$ points  orbit
configuration space of the $\Z_7$-space  $S^3 \setminus \Z_7$.
The forgetful map picking the first coordinate defines a bundle $Y_{2} \to
S^3 \setminus \Z_7$ having as fiber $S^3$ with 14 points removed.
By iterating this procedure we find a tower of fibrations
expressing $Y_{n-1}$ as twisted product, up to homotopy, of
the wedges of spheres $\vee_6 S^2, \vee_{13} S^2$, and so on.
The additive homology of $Y_{n-1}$ splits as tensor product of the 
homology of the factors, by the Serre spectral sequence.
In particular there is a map $\vee_{(n-1)(7n-2)/2} S^2
\to Y_{n-1}$ inducing isomorphism on $H_2$.
The product map $S^3 \times \vee_{(n-1)(7n-2)/2} S^2
\to \widetilde{F}_n(L_{7,1})$ induces isomorphism
on the cohomology groups $H^2, H^3, H^5$. Thus  
all Massey products on elements of $H^2(\widetilde{F}_n(L_{7,1}))$ must vanish.

The unordered configuration space $C_n(L_{7,j})=F_n(L_{7,j})/\Sigma_n$ has
as fundamental group the wreath product $\Sigma_n \wr \Z_7$ 
and has the same universal cover as the ordered configuration space.
It follows that also all unordered configuration spaces are not homotopy
invariant.

Our approach shows that infinite other pairs of homotopic lens spaces
have non homotopic configuration spaces.
It might be interesting to study whether the homotopy type of 
configuration spaces distinguishes up to homeomorphism all lens spaces.

\end{document}